\documentclass[12pt]{article}
\usepackage{amssymb}
\usepackage{amsmath}
\usepackage{amsfonts}
\usepackage{amsthm}
\newtheorem{theorem}{Theorem}
\newtheorem{proposition}{Proposition}
\newtheorem{definition}{Definition}
\newtheorem{lemma}{Lemma}
\newtheorem{corollary}{Corollary}
\theoremstyle{remark}
\newtheorem{remark}{\bf Remark}

\DeclareMathOperator{\Rep}{Rep}
\DeclareMathOperator{\Ob}{Ob}
\DeclareMathOperator{\Mor}{Mor}

\def\address{The University of Iowa, Iowa City, Iowa, 52242-1419, USA\\
Kyiv Taras Shevchenko University,
Cybernetics Department, Volodymyrska, 64, Kyiv, 01033, Ukraine\\
Institute of Mathematics NAS of Ukraine, Tereshchenkivska, 3,
01601, Kyiv, Ukraine}
\def\email{jorgen@math.uiowa.edu\\prosk@unicyb.kiev.ua\\ yurii\_sam@imath.kiev.ua}

\begin{document}
\title{On $C^*$-algebras generated by pairs of $q$-commuting isometries}
\author{Palle E.T. J{\o}rgensen\thanks{This material is based upon work supported by the U.S.
National Science Foundation under Grant No.\ DMS-0139473 (FRG).}, Daniil P. Proskurin and
Yuri\u\i{} S. Samo\u\i{}lenko}
\date{}
\maketitle
\begin{abstract}
We consider the $C^*$-algebras $\mathcal{O}_2^q$ and
$\mathcal{A}_2^q$ generated,
respectively,
by isometries $s_1,s_2$ satisfying
the relation $s_1^*s_2=qs_2s_1^*$ with $\left\vert q\right\vert<1$ (the deformed Cuntz relation),
and by isometries $s_1,s_2$ satisfying
the relation $s_2s_1=qs_1s_2$ with $\left\vert q\right\vert=1$.
We show
that $\mathcal{O}_2^q$ is isomorphic to the Cuntz-Toeplitz
$C^*$-algebra $\mathcal{O}_2^0$ for any $\left\vert q\right\vert <1$.
We further prove that $\mathcal{A}_2^{q_1}\simeq\mathcal{A}_2^{q_2}$ if and only if
either $q_1=q_2$ or $q_1={\overline{q}}_2$.
In the second part of our paper, we discuss
the complexity of the
representation theory of $\mathcal{A}_2^q$.
We show that $\mathcal{A}_2^q$ is $*$-wild
for any $q$ in the circle $\left\vert q\right\vert=1$, and hence
that $\mathcal{A}_2^q$ is not nuclear
for any $q$ in the circle.
\end{abstract}
\section*{Introduction}
The general question of deformations of algebras has received
considerable attention in mathematical physics, and in operator algebra
theory. The motivation derives from commutation relations of quantum
mechanics, see e.g., \cite{osam}, and two distinct cases have received special
attention.
The first case
 (a)
is that of
algebras built on rotation models, often called rotation
algebras; see \cite{riff}. The algebras are generated by a finite set of unitaries,
and a multiplicative commutation relation. These generalized
rotation algebras are $C^*$-algebras, and they are labeled by a rotation
number, or a rotation matrix, in any case labels depending on continuous
parameters. Typically, different parameters yield non-isomorphic $C^*$-algebras.
For example, different irrational rotation numbers yield non-isomorphic
simple $C^*$-algebras. The second case
 (b)
is that of
the additive version of the classical quantum commutation
relations.
This case invites similar deformation questions,
but the answers stand in sharp contrast:
in case (a) we have distinct isomorphism classes, while in case (b)
we have isomorphic $C^*$-algebras for the parameter $q$ in an open set.
In \cite{jsw1}
and \cite{jsw} the authors showed that the Cuntz algebras, see
\cite{cun1}, fit this picture, and that the cases of the canonical
commutation relations (CCR) and the canonical anticommutation
relations (CAR) arise from this viewpoint as limiting cases. It is
shown in \cite{jsw1} that for each $n$ there is a parameter
interval $J$, independent of $n$, containing $0$  but
smaller than  $-1 < q < 1$,
for which a family of
$C^*$-algebras $\mathcal{E}_q$
can be constructed from $n$ generators, so that  for $q =
0$
$\mathcal{E}_{q}$ is $\mathcal{O}_n^0$, the Cuntz-Toeplitz algebra,
and the family $\{ \mathcal{E}_q \mid q \in J \}$
represents only one isomorphism class.
The limiting cases $q=-1$ and
$q=1$ correspond to algebras of the CCRs and the CARs, and the
question was raised in \cite{jsw1} whether the $C^*$-algebras
$\mathcal{E}_q$ are in fact isomorphic for all $q$ in the
``natural'' interval, $-1 < q < 1$.
This has since become a conjecture of some standing,
see, e.g., \cite{mar}.
We show that there is a closely related deformation
family,
and we resolve
the question for this deformation
in the affirmative, in the case $n=2$.

The relationship between the separate deformation systems is
clarified as follows:
The deformation of CCR named $q_{ij}$-CCR
was first constructed by M. Bozejko and R. Speicher \cite{BoSp91}.
It is the $C^*$-algebra generated by $a_i,$ $a_i^*$, $i=1,\dots,d$, satifying the relations
\[
a_i^*a_i=1+q_ia_ia_i^*,\quad a_i^*a_j=q_{ij}a_ja_i^*,\; i<j,\;  q_i\in (-1,1),\; \left| q_{ij}\right|\le 1.
\]
The general conjecture of Jorgensen, Schmitt and Werner \cite{jsw1}
states that for any $q_i\in (-1,1)$, $\left| q_{ij}\right|<1$, 
$i,j=1,\dots,d$, the $C^*$-algebra generated by these relations is isomorphic to the Cuntz-Toeplitz algebra. 

Some limiting cases of the parameters
were recently considered by K. Yuschenko,
who showed that the $C^*$-algebra generated 
by elements satisfying the relations
\[
a_i^*a_i=1+q_ia_ia_i^*,\; q_i\in (0,1),\quad a_i^*a_j=0,\; i\neq j,
\]
is isomorphic to the algebra $\mathcal{O}_n^0$. We consider another limit situation,
putting $q_i=0$. Then we get the algebras 
$\mathcal{O}_n^{q_{ij}}$ generated by isometries $s_i$ satisfying relations of the form
\begin{equation}\label{qizerolim}
s_i^*s_j=q_{ij}s_js_i^*,\; \left| q_{ij}\right|<1,\; i<j.
\end{equation}
The conjecture, in this case,
says that for any $\left| q_{ij}\right|<1$ one should get the Cuntz-Toeplitz algebra.
In this paper we give the affirmative answer for the
case of two generators. We prove isomorphism of the two
$C^*$-algebras $\mathcal{O}_2^q$ and $\mathcal{O}_2^0$ for any $\left| q\right|<1$.

A special case of (\ref{qizerolim}) is the case of
the algebras
generated by $q$-commuting isometries. It was shown in
\cite{dpl} that the $C^*$-algebras generated by the generalized
quonic relations, see \cite{mar}, can be generated by isometries
satisfying relations of the form
\begin{equation}\label{intr}
s_i^*s_j=qs_js_i^*,\quad 1\le i<j\le n,\qquad \left\vert q\right\vert=1.
\end{equation}
We consider here the situation with $\left\vert q\right\vert<1$
which is also the deformation of the Cuntz-Toeplitz algebra. Let
us denote the $C^*$-algebras generated by isometries satisfying
(\ref{intr}) by $\mathcal{O}_n^q$. Then using the general methods
developed in \cite{jsw1} one can find some $0<\varepsilon<1$ such
that for any $\left\vert q\right\vert<\varepsilon$ one has
$\mathcal{O}_n^q\simeq\mathcal{O}_n^0$. And again one has the
question whether the isomorphism exists for any $q$ with $\left\vert
q\right\vert<1$.

In this paper we show that this is the case for $n =2$. The
methods which we use are different from the general case of $n>
2$, and in any case are of independent interest. The algebras
$\mathcal{O}_2^q$ have natural Fock representations, acting on an
infinite ``particle'' Hilbert space, constructed by use of the
Fock tensor functor. An important issue, which we resolve, is the
faithfulness of the Fock representation.

In Section \ref{SecCsaO1}, we study the $C^*$-algebra
$\mathcal{O}_2^q$ generated by isometries $s_1,s_2$ satisfying the
deformed Cuntz relation
\[
s_1^*s_2=qs_2s_1^*,\quad \left\vert q\right\vert<1.
\]
Indeed, for $q=0$, we get the Cuntz-Toeplitz $C^*$-algebra
$\mathcal{O}_2^0$, see \cite{cun1}, generated by isometries
$t_1,t_2$ satisfying the relation of orthogonality $t_1^*t_2=0$.

We prove, see Sec.\ \ref{SubsecIsom1}, that $\mathcal{O}_2^q\simeq
\mathcal{O}_2^0$,
$C^*$-isomorphism,
for any $q$, $\left\vert q\right\vert<1$. In Sec.\ \ref{SubsecRepr2} we
construct special representations of $\mathcal{O}_2^q$ and discuss
on an informal level the constructions presented in Sec.\ \ref{SubsecIsom1}.

The situation with $\left\vert q\right\vert=1$ is quite different, see
\cite{kab,dpl,rmp}. Almost all $\mathcal{O}_2^q$ with $\left\vert
q\right\vert=1$ are non-isomorphic. To be more precise,
$\mathcal{O}_2^{q_1}\simeq\mathcal{O}_2^{q_2}$ if and only if the
corresponding non-commutative tori $\mathcal{A}_{q_j}$, $j=1,2$,
are isomorphic. Recall that the non-commutative torus
$\mathcal{A}_q$, $\left\vert q\right\vert=1$, is the $C^*$-algebra generated by
a pair of unitaries $u_1,u_2$ satisfying the relation
\[
u_2u_1=qu_1u_2 ;
\]
see, for example, \cite{riff}. Put $q_j=e^{2 i \pi\theta_j}$,
$\theta_j\in\mathbb{R}$, $j=1,2$; then
$\mathcal{A}_{q_1}\simeq\mathcal{A}_{q_2}$ if and only if
$\theta_2=\pm\theta_1 \pmod{\mathbb{Z}}$.

It is easy to see that the relation $s_1^*s_2= qs_2s_1^*$, $\left\vert
q\right\vert=1$, implies the relation $s_2s_1= q s_1s_2$. Indeed, for
$a=s_2s_1-qs_1s_2$ one has $a^*a=0$. However, the converse is not
true. In particular, $\mathcal{O}_2^q$ is nuclear for any $q$,
$\left\vert q\right\vert= 1$, but one of the results of Section \ref{SecCsaA2} implies that
the $C^*$-algebra generated by isometries $s_1,s_2$ satisfying
$s_2s_1=qs_1s_2$ is not nuclear.

In Section \ref{SecCsaA2}, we consider the $C^*$-algebra $\mathcal{A}_2^q$
generated by isometries $s_1,s_2$ satisfying the relation
\[
s_2s_1=qs_1s_2,\quad \left\vert q\right\vert=1 .
\]
The relation $s_2s_1=q s_1s_2$ implies that $\left\vert q\right\vert=1$.
Indeed, let the isometries be realized on a Hilbert space
$\mathcal{H}$: then for any non-zero $x\in\mathcal{H}$ one has
\[
\left\Vert x\right\Vert=\left\Vert s_2s_1 x\right\Vert= \left\vert q\right\vert \left\Vert s_1s_2 x\right\Vert
=\left\vert q\right\vert\left\Vert x\right\Vert.
\]

It is proved, see Sec.\ \ref{SubsecDesc1}, that
$\mathcal{A}_2^{q_1}\simeq\mathcal{A}_2^{q_2}$ if and only if either
$q_1=q_2$ or $q_1={\overline{q}}_2$. Further, in Sec.\ \ref{SubsecAwild} we
show that the representation theory of $\mathcal{A}_2^q$ is
extremely complicated, even in the case $ q=1$~! More precisely,
the problem of the classification of the irreducible
representations of $\mathcal{A}_2^q$ contains as a sub-problem
the description of the irreducible representations of
$C^*(\mathcal{F}_2)$, where $C^*(\mathcal{F}_2)$ is the group
$C^*$-algebra of the free group with two generators. The
classification of irreducible representations of
$C^*(\mathcal{F}_2)$, or equivalently the description of all
irreducible pairs of unitary operators, is the standard ``$*$-wild''
problem; see \cite{osam} for a detailed discussion of the
complexity of the representation theory of $C^*$-algebras.

The properties of pairs of commuting proper isometries were originally studied in
\cite{bcl}, where a construction demonstrating the
$*$-wildness of the $C^*$-algebra $\mathcal{A}_2^1$ was presented.
Note that the construction which we use to prove that
$\mathcal{A}_2^q$ is $*$-wild is not a generalization of the one
presented in \cite{bcl}.

\section{\label{SecCsaO1}The $C^*$-algebra $\mathcal{O}_2^q$}
\subsection{\label{SubsecIsom1}The isomorphism $\mathcal{O}_2^q\simeq\mathcal{O}_2^0$}
In this part we show that for any $q$, $\left\vert q\right\vert<1$, there is an
isomorphism $\mathcal{O}_2^q\simeq\mathcal{O}_2^0$,
i.e., isomorphism of $C^*$-algebras. It is a
special case of the hypothesis of P.E.T. J{\o}rgensen, L.M.
Schmitt and R.F. Werner presented in \cite{jsw1,jsw}. Namely, it
was shown that the $C^*$-algebras defined by generators
$a_i,a_i^*$, $i=1,\ldots,d$, satisfying the relations
\[
a_i^*a_j=\delta_{ij}1+\sum_{k,l=1}^d T_{ij}^{kl}a_la_k^*,\quad
T_{ij}^{kl}=\overline{T}_{ji}^{lk}\in\mathbb{C} ,
\]
for sufficiently small absolute values of the coefficients, are
isomorphic to the Cuntz-Toeplitz $C^*$-algebra $\mathcal{O}_d$,
see \cite{cun1}, generated by isometries $t_i$, $i=1,\ldots,d$,
satisfying the relations
\[
t_i^*t_j=0,\quad i\ne j .
\]
To be more precise, the norm bound $\left\Vert T\right\Vert<\sqrt{2}-1$ gives
a sufficient condition for this isomorphism. Here $T$ is the
self-adjoint operator acting on $\mathcal{H}^{\otimes 2}$,
$\mathcal{H}=\bigl<e_1,\ldots,e_d\bigr>$, defined by the action on
the basis as follows:
\[
T e_{k}\otimes e_{l} = \sum_{i,j} T_{ik}^{lj}e_{i}\otimes e_{j}.
\]
It was conjectured in \cite{jsw} that the result is correct for
$\left\Vert T\right\Vert<1$; see also \cite{dn}.

Below we prove that, for all $q$ in the open interval $-1<q<1$, the $C^*$-algebra
$\mathcal{O}_2^q$ can be generated by generators of
$\mathcal{O}_2^0$, and vice versa.
\begin{remark}
All of the arguments presented below carry over to the case of any $q$
in the complex disk $\left\vert q\right\vert<1$.
\end{remark}
We prove our result by several lemmas.
\begin{lemma}
\label{Lem1}
Let $s_1,s_2$ be isometries satisfying the relation
$s_1^*s_2=qs_2s_1^*$ with $-1<q<1$. Construct the elements
\begin{align*}
t_1(s_1,s_2)=t_1&:=s_1,\\
t_2(s_1,s_2)=t_2&:=(1-s_1s_1^*)s_2(1-s_1s_1^*+(1-q^2)^{-\frac{1}{2}}s_1s_1^*).
\end{align*}
Then $t_i^*t_i=1$, $i=1,2$, and $t_1^*t_2=0$.
\end{lemma}
\begin{proof}
The relation $t_1^*t_2=0$ follows from the relation $s_1^*(1-s_1s_1^*)=0$.
Let us verify that $t_2^*t_2=1$.
\begin{align*}
t_2^*t_2&=(1-s_1s_1^*+(1-q^2)^{-\frac{1}{2}}s_1s_1^*)s_2^*(1-s_1s_1^*)s_2
(1-s_1s_1^*+(1-q^2)^{-\frac{1}{2}}s_1s_1^*)\\ &=
(1-s_1s_1^*+(1-q^2)^{-\frac{1}{2}}s_1s_1^*)(1-s_1s_1^*+(1-q^2)^{-\frac{1}{2}}s_1s_1^*)\\
&\qquad -(1-s_1s_1^*+(1-q^2)^{-\frac{1}{2}}s_1s_1^*)s_2^*s_1s_1^*s_2
(1-s_1s_1^*+(1-q^2)^{-\frac{1}{2}}s_1s_1^*)\\ &=
1+\frac{q^2}{1-q^2}s_1s_1^*\\
&\qquad -
q^2(1-s_1s_1^*+(1-q^2)^{-\frac{1}{2}}s_1s_1^*)s_1s_1^*
(1-s_1s_1^*+(1-q^2)^{-\frac{1}{2}}s_1s_1^*)\\ &=
1+\frac{q^2}{1-q^2}s_1s_1^*-\frac{q^2}{1-q^2}s_1s_1^*=1,
\end{align*}
where we used the relation
$s_2^*s_1s_1^*s_2=q^2s_1s_2^*s_2s_1^*=q^2s_1s_1^*$.
\end{proof}
In the following lemma we present the converse construction.
\begin{lemma}
\label{Lem2}
Let $t_1,t_2$ be isometries satisfying $t_1^*t_2=0$. Construct the
operators
\[
s_1(t_1,t_2)=s_1:=t_1,\quad
\widetilde{t}_2:=t_2(1-t_1t_1^*+(1-q^2)^{\frac{1}{2}}t_1t_1^*)
\]
and put
\[
s_2(t_1,t_2)=s_2:=\sum_{i=0}^{\infty}q^i
t_1^i\widetilde{t_2}(t_1^*)^i,\quad t_1^0:=1.
\]
Then $s_i^*s_i=1$, $i=1,2$, and $s_1^*s_2=qs_2s_1^*$.
\end{lemma}
\begin{proof}
Evidently for $-1<q<1$ the series converges with respect to norm.
Let us show that $s_1^*s_2=qs_2s_1^*$. We will use the obvious
relation $t_1^*\widetilde{t}_2=0$.
\begin{align*}
s_1^*s_2& =t_1^*\sum_{i=0}^{\infty}q^i
t_1^i\widetilde{t_2}(t_1^*)^i= \sum_{i=1}^{\infty}q^i
t_1^{i-1}\widetilde{t_2}(t_1^*)^i\\
&=q\left(\sum_{i=0}^{\infty}q^i
t_1^i\widetilde{t_2}(t_1^*)^i\right)t_1^*=qs_2s_1^*.
\end{align*}
To show that $s_2^*s_2=1$ we note firstly that
\[
\widetilde{t}_2^*\widetilde{t}_2=1-q^2t_1t_1^*
\]
and
\[
t_1^i\widetilde{t_2^*}(t_1^*)^it_1^l\widetilde{t}_2(t_1^*)^l=\left\{
\begin{aligned}
{} & t_1^i\widetilde{t_2^*}(t_1^*)^{i-l}\widetilde{t}_2(t_1^*)^l& &=0,& &
i>l , \\
{} & t_1^i\widetilde{t_2^*}t_1^{l-i}\widetilde{t}_2(t_1^*)^l& &=0,& &
l>i , \\
{} & \rlap{$\displaystyle t_1^i(1-q^2t_1t_1^*)(t_1^i)^*,$} & & & & l=i .
\end{aligned}
\right.
\]
Then
\begin{align*}
s_2^*s_2&=\bigl(\sum_{i=0}^{\infty}q^i
t_1^i\widetilde{t_2}^*(t_1^*)^i\bigr)
\bigl(\sum_{l=0}^{\infty}q^lt_1^l\widetilde{t_2}(t_1^*)^l\bigr)\\&=
\sum_{i=0}^{\infty}q^i t_1^i(1-q^2t_1t_1^*)(t_1^*)^i=
\sum_{i=0}^{\infty}\bigl(q^{2i}
t_1^i(t_1^*)^i-q^{2i+2}t_1^{i+1}(t_1^*)^{i+1}\bigr)=1.
\end{align*}
\end{proof}
In the following lemma we show that starting from generators of
$\mathcal{O}_2^q$ and applying consecutively the constructions
presented in Lemmas \ref{Lem1} and \ref{Lem2}, we get the starting elements.
\begin{lemma}
\label{Lem3}
\[
s_i(t_1(s_1,s_2),t_2(s_1,s_2))=s_i,\quad i=1,2.
\]
\end{lemma}
\begin{proof}
For
$t_2=(1-s_1s_1^*)s_2(1-s_1s_1^*+(1-q^2)^{-\frac{1}{2}}s_1s_1^*)$,
one has
\begin{align*}
\widetilde{t}_2&=(1-s_1s_1^*)s_2(1-s_1s_1^*+(1-q^2)^{-\frac{1}{2}}s_1s_1^*)
(1-s_1s_1^*+(1-q^2)^{\frac{1}{2}}s_1s_1^*)\\&=(1-s_1s_1^*)s_2.
\end{align*}
Note that
$s_1^i(1-s_1s_1^*)s_2(s_1^*)^i=s_1^is_2(s_1^*)^i-qs_1^{i+1}s_2(s_1^*)^{i+1}$.
Then we get
\[
s_2(t_1,t_2)=\sum_{i=0}^{\infty}q^i s_1^i\widetilde{t}_2(s_1^*)^i=
\sum_{i=0}^{\infty}\left(q^is_1^is_2(s_1^*)^i-q^{i+1}s_1^{i+1}s_2(s_1^*)^{i+1}\right)=s_2.
\]
\end{proof}
In fact we have proved our result.
\begin{theorem}\label{thst}
The isomorphism
$\mathcal{O}_2^q\simeq\mathcal{O}_2^0$
holds for all $q$ in the open interval $-1<q<1$.
\end{theorem}
\begin{proof}
The isomorphism
$\phi\colon\mathcal{O}_2^0\rightarrow\mathcal{O}_2^q$ is defined
by
\[
\phi(t_1)=s_1,\ \phi(t_2)=t_2(s_1,s_2).
\]
It follows from Lemma \ref{Lem3} that the inverse homomorphism is given by
the formulas
\[
\psi(s_1)=t_1,\quad \psi(s_2)=s_2(t_1,t_2).
\]
\end{proof}
\subsection{\label{SubsecRepr2}Representations of $\mathcal{O}_2^q$}
In this part we recall the notion of the Fock representation of
our $q$-relations and construct a special class of
representations of $\mathcal{O}_2^q$ which includes the Fock one.
Using these special representations we discuss informally the
results of Sec.\ \ref{SubsecIsom1}.

Recall that the Fock representation of $\mathcal{O}_2^0$
($\mathcal{O}_2^q$) is the unique irreducible representation
defined by the vacuum vector $\Omega$ with the property
$t_i^*\Omega=0$ ($s_i^*\Omega=0$), $i=1,2$. It follows from the
main result of \cite{jps} that the Fock representation of
$\mathcal{O}_2^q$ is positive, i.e., it is a $*$-representation on
the Hilbert space, and faithful at the algebraic level. The later
means that the Fock representation of the $*$-algebra generated by
the basic relations of $\mathcal{O}_2^q$ has trivial kernel.
Our stability theorem implies that the same is at the
$C^*$-algebra level.
\begin{proposition}
\label{Prop1}
The Fock representation of $\mathcal{O}_2^q$ is faithful.
\end{proposition}
\begin{proof}
It is known that the Fock representation of the
Cuntz-Toeplitz algebra is faithful. This follows from the fact that
any irreducible representation of the Cuntz-Toeplitz algebra is
either the Fock representation or a representation of the Cuntz
algebra (see \cite{cun1}). Hence we have to show that the Fock
representation of $\mathcal{O}_2^0$ is the Fock representation of
$\mathcal{O}_2^q$, i.e.,  if $t_i^*\Omega=0$, $i=1,2$, then
$s_i^*\Omega=0$, $i=1,2$. Since $s_1=t_1$ we have to verify only
that $s_2^*\Omega=0$. Indeed, the conditions $t_i^*\Omega=0$, $i=1,2$, imply that
$\widetilde{t}_2^*\Omega=0$ and $s_2(t_1,t_2)^*\Omega=0$.
\end{proof}
To clarify the nature of the constructions presented in Sec.\ \ref{SubsecIsom1}, we
consider a special class of representations of $\mathcal{O}_2^q$. Put
$s_1$ to be a multiple of the unilateral shift, i.e., suppose
that the representation space is
$l_2(\mathbb{N})\otimes\mathcal{K}$ and that
\[
s_1=S\otimes\mathbf{1}=\left(
\begin{array}{cccc}
0&0&0&\cdots\\ \mathbf{1}&0&0&\cdots\\ 0&\mathbf{1}&0&\cdots\\
&\ddots&\ddots&\ddots
\end{array}
\right) .
\]
Then it is easy to see that the relation $s_1^*s_2=qs_2s_1^*$
implies that $s_2$ has the matrix form
\begin{equation}\label{qcun}
s_2=\left(
\begin{array}{ccccc}
u_1&\sqrt{1-q^2}u_2&\sqrt{1-q^2}u_3& \sqrt{1-q^2}u_4&\cdots
\\
0&qu_1&q\sqrt{1-q^2}u_2&q\sqrt{1-q^2}u_3&\cdots
\\
0&0&q^2u_1&q^2\sqrt{1-q^2}u_2&\cdots\\
&\ddots&\ddots&\ddots&\ddots\\
\end{array}
\right) ,
\end{equation}
where the elements $u_i$, $i\in\mathbb{N}$, are generators of the
Cuntz algebra $\mathcal{O}_{\infty}$, i.e., satisfy the relations
\[
u_i^*u_i=1,\quad u_i^*u_j=0,\ i\ne j .
\]
Analogously, for generators of $\mathcal{O}_2^0$ we have
\begin{equation}\label{cun}
t_1=\left(
\begin{array}{cccc}
0&0&0&\cdots\\ \mathbf{1}&0&0&\cdots\\ 0&\mathbf{1}&0&\cdots\\
&\ddots&\ddots&\ddots
\end{array}
\right),\quad
t_2=\left(
\begin{array}{ccccc}
u_1&u_2&u_3&u_4&\cdots
\\
0&0&0&0&\cdots
\\
0&0&0&0&\cdots\\ &\ddots&\ddots&\ddots&\ddots\\
\end{array}
\right),
\end{equation}
where the elements $u_i$, $i\in\mathbb{N}$, satisfy the same
relations as in the $q$-deformed case. Then the constructions
presented above are just the elementary matrix transformations
reducing the matrix of the form (\ref{qcun}) to the matrix of the
form (\ref{cun}) and vice versa.

Note that the presented motivations cannot be treated as the
correct proof of Theorem \ref{thst} since \emph{a priori} it is not known
whether or not the constructed representation of $\mathcal{O}_2^q$ is
faithful.

In fact, the construction (\ref{cun}) determines the  functor
\[
F\colon \Rep \mathcal{O}_{\infty}\rightarrow \Rep \mathcal{O}_2^0,
\]
see the next section, defined as follows. For any $\pi\in \Ob(\Rep
\mathcal{O}_{\infty})$, the representation $F(\pi)\in \Ob(\Rep
\mathcal{O}_2^0)$ is constructed by formulas (\ref{cun}) and for
any $\Lambda\in \Mor(\pi_1,\pi_2)$ one has
$F(\Lambda)=\mathbf{1}\otimes\Lambda$. It can be verified that the
functor $F$ is full and faithful. In particular the representation
$F(\pi)$ is irreducible if and only if the representation $\pi$ is
irreducible, and $F(\pi_1)$ is unitarily equivalent to $F(\pi_2)$ if and only if
$\pi_1$ and $\pi_2$ are unitarily equivalent.
\begin{remark}
It is easy to see that the Fock representation of
$\mathcal{O}_2^0$ corresponds to the Fock representation of
$\mathcal{O}_{\infty}$.
\end{remark}
\section{\label{SecCsaA2}The $C^*$-algebra $\mathcal{A}_2^q$}
\subsection{\label{SubsecDesc1}Description of the isomorphism classes}
In this part we show that
$\mathcal{A}_2^{q_1}\simeq\mathcal{A}_2^{q_2}$ if and only if either
$q_1=q_2$ or $q_1=\overline{q}_2$. In  the following it will be
convenient for us to put $q$ in the form $q=e^{2 i \pi\theta}$.
\begin{proposition}
\label{Prop2}
Let $\theta$ be irrational, $q=e^{2 i \pi\theta}$. Then there
exists a unique normalized tracial state $\widetilde{\tau}$ on
$\mathcal{A}_2^q$, and $\mathcal{M}=\{a \mid
\widetilde{\tau}(a^*a)=0\}$ is the largest two-sided closed ideal
in $\mathcal{A}_2^q$.
\end{proposition}
\begin{proof}
Recall that by $\mathcal{A}_q$ we denote the non-commutative torus
corresponding to $q$. Let $\mathcal{M}$ be the closed two-sided
ideal generated by the projections $1-s_js_j^*$, $j=1,2$. Then one
has the canonical homomorphism
\[
\varphi\colon\mathcal{A}_2^q\rightarrow
\mathcal{A}_2^q/\mathcal{M}=\mathcal{A}_q.
\]
Since $\mathcal{A}_q$ with irrational $\theta$ is simple, see
\cite{riff}, $\mathcal{M}$ is the largest ideal in
$\mathcal{A}_2^q$. Put $\widetilde{\tau}=\tau\circ\varphi$, where
$\tau$ is the unique normalized tracial state on $\mathcal{A}_q$,
see \cite{riff}. Since
\[
\mathcal{J}=\{a
\mid
\widetilde{\tau}(a^*a)=0\}
\]
is the two-sided closed ideal and
$\widetilde{\tau}(1-s_js_j^*)=0$, $j=1,2$, one has
$\mathcal{M}\subset\mathcal{J}$. Hence $\mathcal{M}=\mathcal{J}$.

Let $\widetilde{\sigma}$ be a normalized tracial state on
$\mathcal{A}_2^q$; then as above,
\[
\widetilde{\mathcal{M}}:=\{a
\mid
\widetilde{\sigma}(a^*a)=0\}=\mathcal{M},
\]
and
in particular $\widetilde{\sigma}(a)=0$ for any $a\in\mathcal{M}$.
Then one can define the tracial state $\sigma$ on $\mathcal{A}_q$
by the rule $\sigma(b)=\widetilde{\sigma}(a)$, $\varphi(a)=b$.
Evidently $\widetilde{\sigma}=\sigma\circ\varphi$. By the
uniqueness of the trace on $\mathcal{A}_q$, one has $\tau=\sigma$
and $\widetilde{\sigma}=\widetilde{\tau}$.
\end{proof}
As a corollary, we have the isomorphism of
$\mathcal{A}_2^{q_j}$ with irrational $\theta_j$, $j=1,2$.
\begin{proposition}
\label{Prop3}
Consider $\mathcal{A}_2^{q_j}$, $q_j=e^{2 i \pi\theta_j}$, with
irrational $\theta_j$, $j=1,2$. Then
$\mathcal{A}_2^{q_1}\simeq\mathcal{A}_2^{q_2}$ if and only if
 $\theta_1=\pm\theta_2 \pmod{\mathbb{Z}}$.
\end{proposition}
\begin{proof}
We prove that the isomorphism
$\mathcal{A}_2^{q_1}\simeq\mathcal{A}_2^{q_2}$ implies
$\theta_1=\pm\theta_2 \pmod{\mathbb{Z}}$. The converse
implication is trivial.

Let $\psi\colon\mathcal{A}_2^{q_1}\rightarrow\mathcal{A}_2^{q_2}$
be an isomorphism and $\widetilde{\tau}_j$, $j=1,2$, be normalized
traces on $\mathcal{A}_2^{q_j}$, $j=1,2$. Denote by
$\mathcal{M}_j\subset\mathcal{A}_2^{q_j}$, $j=1,2$, the largest
ideals introduced in the proposition above.

Consider the normalized tracial state
$\widehat{\tau}_1=\widetilde{\tau}_2\circ\psi$. Then
$\widehat{\tau}_1=\widetilde{\tau}_1$, and for any
$a\in\mathcal{M}_1$ one has
\[
\widehat{\tau}_1(a^*a)=\widetilde{\tau}_2(\psi(a^*a))=0,
\]
i.e., $\psi(a)\in\mathcal{M}_2$ and
$\psi(\mathcal{M}_1)\subset\mathcal{M}_2$. Analogously
$\psi^{-1}(\mathcal{M}_2)\subset\mathcal{M}_1$. Hence $\psi$
induces the isomorphism
\[
\psi\colon\mathcal{M}_1\rightarrow\mathcal{M}_2.
\]
Denote by $\phi_2$ the canonical homomorphism
\[
\phi_2\colon\mathcal{A}_2^{q_2}\rightarrow\mathcal{A}_2^{q_2}/\mathcal{M}_2=
\mathcal{A}_{q_2}\,;
\]
then $\ker \phi_2\circ\psi=\mathcal{M}_1$, hence
$\mathcal{A}_{q_1}\simeq\mathcal{A}_{q_2}$ and $\theta_1=\pm
\theta_2 \pmod{\mathbb{Z}}$.
\end{proof}
Let us now consider rational $\theta_j$, $j=1,2$.
\begin{proposition}
\label{Prop4}
Let $\theta_j\in\mathbb{Q}$, $j=1,2$. Then
$\mathcal{A}_2^{q_1}\simeq \mathcal{A}_2^{q_2}$ if and only if
 $\theta_1=\pm\theta_2 \pmod{\mathbb{Z}}$.
\end{proposition}
\begin{proof}
As in the irrational case, we prove that the existence of an isomorphism
\[
\psi\colon\mathcal{A}_2^{q_1}\rightarrow\mathcal{A}_2^{q_2}
\]
implies isomorphism of $\mathcal{A}_{q_1}$ and $
\mathcal{A}_{q_2}$.

Let $\mathcal{M}_j\subset\mathcal{A}_2^{q_j}$, $j=1,2$, be the
ideals generated by the projections $1-s_is_i^*$, $i=1,2$. We show
that $\psi(\mathcal{M}_1)\subset\mathcal{M}_2$. Indeed, consider
the canonical homomorphism
\[
\varphi_2\colon\mathcal{A}_2^{q_2}\rightarrow\mathcal{A}_{q_2}
\]
and the composite homomorphism
\[
\varphi_2\circ\psi\colon\mathcal{A}_2^{q_1}\rightarrow
\mathcal{A}_{q_2}.
\]
Since any irreducible representation of $\mathcal{A}_{q}$ with
rational $\theta$ is finite-dimensional, $\mathcal{A}_{q}$
does not contain any non-unitary isometry. Hence
$\varphi_2\circ\psi(1-s_is_i^*)=0$ and
$\psi(1-s_is_i^*)\in\mathcal{M}_2$, $i=1,2$. So
$\psi(\mathcal{M}_1)\subset\mathcal{M}_2$. Analogously,
$\psi^{-1}(\mathcal{M}_2)\subset\mathcal{M}_1$. As in the proof of
the previous proposition, one has
$\mathcal{A}_{q_1}\simeq\mathcal{A}_{q_2}$. It remains only to recall
that the rational tori $\mathcal{A}_{q_j}$, $j=1,2$, are
isomorphic if and only if $\theta_1=\pm\theta_2 \pmod{\mathbb{Z}}$.
\end{proof}
\subsection{\label{SubsecAwild}The $C^*$-algebra $\mathcal{A}_2^q$ is $*$-wild}
In this part we discuss the complexity of the representation
theory of $\mathcal{A}_2^q$.

To  compare $C^*$-algebras according to the complexity of
their categories of representations, we use the relation of
majorization. Note that we modify the definition of majorization
given in \cite{osam} to make it less restrictive.

Recall that the category of representations of a certain
$C^*$-algebra $\mathcal{A}$, denoted by $\Rep \mathcal{A}$, has
the representations of $\mathcal{A}$ as its objects and the
intertwining operators as its morphisms.
\begin{definition}
\label{major}
\upshape
We say that a $C^*$-algebra $\mathcal{A}$ is \emph{majorized} by a
$C^*$-algebra $\mathcal{B}$, $\mathcal{A}\prec\mathcal{B}$, if
there exist a homomorphism
\[
\varphi\colon\mathcal{B}\rightarrow\mathcal{A}\otimes\mathcal{C},
\]
where $\mathcal{C}$ is a nuclear $C^*$-algebra, and an
irreducible representation
\[
\widetilde{\pi}\colon\mathcal{C}\rightarrow B(\mathcal{H})
\]
such that the functor
\[
\mathcal{F}_{\varphi}\colon \Rep \mathcal{A}\rightarrow \Rep
\mathcal{B}
\]
defined by
\begin{align*}
\mathcal{F}_{\varphi}(\pi) &
=(\pi\otimes\widetilde{\pi})\circ\varphi,\quad \pi\in \Ob (\Rep
\mathcal{A}),\\ \mathcal{F}_{\varphi}(A) & = A\otimes 1,\quad A\in
\Mor (\pi_1,\pi_2),
\end{align*}
is full and faithful.
\end{definition}

Informally, this definition means that using $\varphi$ one can
construct the representations of $\mathcal{B}$ from the
representations of $\mathcal{A}$ and the representations of
$\mathcal{B}$ are irreducible (unitarily equivalent) if and only if the
corresponding representations of $\mathcal{A}$ are irreducible
(unitarily equivalent).

It is easy to see that majorization is a partial order.
\begin{definition}
\label{starwild}
\upshape
We say that a $C^*$-algebra $\mathcal{A}$ is $*$\emph{-wild} if
\[
C^*(\mathcal{F}_2)\prec\mathcal{A},
\]
where $C^*(\mathcal{F}_2)$ is the group $C^*$-algebra of the free
group with two generators.
\end{definition}
The group $C^*$-algebra $C^*(\mathcal{F}_2)$ is considered as the
standard $*$-wild algebra, since it can be shown that the problem of
the classification of the representations of $C^*(\mathcal{F}_2)$
contains as a sub-problem the classification of the representations
of any finitely generated $C^*$-algebra, in particular
$C^*(\mathcal{F}_n)\prec C^*(\mathcal{F}_2)$; see \cite{osam}.

Obviously a $C^*$-algebra majorizing a $*$-wild $C^*$-algebra is
$*$-wild. 
The $*$-wild $C^*$-algebras have a very complicated category of
representations: in particular, it was noted in \cite{osam} that
$*$-wild algebras are not nuclear. Since our definition of
majorization generalizes the one considered in \cite{osam}, we
present below the proof of this statement.
\begin{proposition}
\label{Prop5}
Let $\mathcal{A}$ be a $*$-wild $C^*$-algebra. Then $\mathcal{A}$ is not nuclear.
\end{proposition}
\begin{proof}
It is sufficient to show that $\mathcal{A}$ has a representation
generating a non-hyperfinite factor, since any
factor-representation of a nuclear $C^*$-algebra is hyperfinite,
by a theorem of Alain Connes.
Since $C^*(\mathcal{F}_2)\prec\mathcal{A}$, one has the
homomorphism
\[
\varphi\colon\mathcal{A}_2^q\rightarrow C^*(\mathcal{F}_2)\otimes
\mathcal{C},
\]
where $\mathcal{C}$ is a nuclear $C^*$-algebra, and the
irreducible representation
\[
\widetilde{\pi}\colon\mathcal{C}\rightarrow B(\mathcal{H}),
\]
as in Definition \ref{major}.

Consider a representation $\pi$ of $C^*(\mathcal{F}_2)$
generating a non-hyperfinite factor. Put
$(\pi\otimes\widetilde{\pi})\circ\varphi:=\pi_1$ and note that
\[
(\pi\otimes\widetilde{\pi})
(C^*(\mathcal{F}_2)\otimes\mathcal{C})^{\prime\prime}=
\pi(C^*(\mathcal{F}_2))^{\prime\prime}\otimes B(\mathcal{H})
\]
is a non-hyperfinite factor also. We use a prime here to denote
the commutant (and double primes for the double commutant).
Since the functor
\[
\mathcal{F}_{\varphi}\colon \Rep C^*(\mathcal{F}_2))\rightarrow
\Rep \mathcal{A}
\]
defined in Definition \ref{major} is full and faithful, one has
\[
\pi_1(\mathcal{A})^{\prime}=(\pi\otimes\widetilde{\pi})
(C^*(\mathcal{F}_2)\otimes\mathcal{C})^{\prime}.
\]
Hence
\[
\pi_1(\mathcal{A})^{\prime\prime}=(\pi\otimes\widetilde{\pi})
(C^*(\mathcal{F}_2)\otimes\mathcal{C})^{\prime\prime}
\]
is a non-hyperfinite factor.
\end{proof}
To prove that $\mathcal{A}_2^q$ is $*$-wild, we need an auxiliary
proposition. Note that below we denote by $\mathcal{O}_2$ the
Cuntz algebra, i.e., we suppose that the generators of
$\mathcal{O}_2^0$ satisfy the additional relation
$t_1t_1^*+t_2t_2^*=1$.

In the following for any $C^*$-algebras $\mathcal{A}_i$, $i=1,2$,
we denote by $\mathcal{A}_1\star\mathcal{A}_2$ their free product,
see \cite{voi}.
\begin{proposition}
\label{main}
The $C^*$-algebra $\mathcal{O}_2\star C([0,1])$ is $*$-wild.
\end{proposition}
\begin{proof}
It was shown in \cite{mtur} that the $C^*$-algebra $C([0,1])\star
C([0,1])$ is $*$-wild. Then to prove our statement it is sufficient
to show that
\[
C([0,1])\star C([0,1])\prec \mathcal{O}_2\star C([0,1]).
\]
In the following we denote by $c$ the standard generator of
$C([0,1])$, $c(x)=x$, for any $x\in [0,1]$, and denote by $c_1,c_2$
the standard free generators of $C([0,1])\star C([0,1])$. Then the
needed majorization is given by the homomorphism
\[
\varphi\colon\mathcal{O}_2\star C([0,1])\rightarrow
\bigl(C([0,1])\star C([0,1])\bigr)\otimes\mathcal{O}_2
\]
defined by
\[
\varphi(x)=1\otimes x,\; x\in\mathcal{O}_2,\quad
\varphi(c)=c_1\otimes t_1t_1^*+c_2\otimes t_2t_2^*
\]
and some irreducible representation
$\widetilde{\pi}\colon\mathcal{O}_2\rightarrow B(\mathcal{H})$.

To prove that the induced functor
\[
\mathcal{F}_{\varphi}\colon \Rep C([0,1])\star C([0,1])\rightarrow
\Rep \mathcal{O}_2\star C([0,1])
\]
is full and faithful, it is sufficient to show that any $\Lambda\in
F(\pi)(\mathcal{O}_2\star C([0,1]))^{\prime}$ has the form
$\Lambda_1\otimes\mathbf{1}$ with $\Lambda_1\in \pi (C([0,1])\star
C([0,1]))^{\prime}$; see Lemma 13 and Remark 49 in \cite{osam}.

Let us put $\mathcal{F}_{\varphi}(\pi):=\pi_1$,
$\widetilde{\pi}(t_i):=T_i$, and $\pi(c_i)=C_i$, $i=1,2$. Then
\[
\pi_1 (c)= C_1\otimes T_1T_1^*+C_2\otimes T_2T_2^*,\quad
\pi_1(t_i)=\mathbf{1}\otimes T_i,\; i=1,2.
\]
Let $\Lambda\in\pi_1(\mathcal{O}_2\star C([0,1]))$. Then it is easy
to see that the relations
\[
\Lambda (\mathbf{1}\otimes T_i)=(\mathbf{1}\otimes T_i)\Lambda,\quad
\Lambda (\mathbf{1}\otimes T_i^*)=(\mathbf{1}\otimes T_i^*)\Lambda
\]
imply that $\Lambda=\Lambda_1\otimes\mathbf{1}$. Further, since
\[
\pi_1(c\, t_it_i^*)=C_i\otimes T_iT_i^*,\quad i=1,2,
\]
one has
\[
\Lambda_1 C_i\otimes T_iT_i^*= C_i\Lambda_1\otimes T_iT_i^*,\quad
i=1,2,
\]
hence $\Lambda_1 C_i=C_i\Lambda_1$, $i=1,2$, and
$\Lambda_1\in\pi(C([0,1])\star C([0,1]))^{\prime}$.
\end{proof}

Now we are ready to prove the main result of this part.
\begin{theorem}
The $C^*$-algebra $\mathcal{A}_2^q$ is $*$-wild.
\end{theorem}
\begin{proof}
We prove that $\mathcal{O}_2\star C([0,1])\prec\mathcal{A}_2^q$.
Let us construct the homomorphism
\[
\varphi\colon\mathcal{A}_2^q\rightarrow (\mathcal{O}_2\star
C([0,1]))\otimes\mathcal{B},
\]
where
\[
\mathcal{B}=C^*(s,u
\mid
s^*s=1,\; u^* u=uu^*=1,\; us=qsu).
\]
This $\mathcal{B}$ is nuclear, since it is the crossed product
$\mathcal{B}=\mathcal{T}(C(\mathbf{T}))\rtimes\mathbb{Z}$, where
$\mathcal{T}(C(\mathbf{T}))$ is the Toeplitz $C^*$-algebra.

Pick a function
$c$ taking values strictly between
$0$ and $1$ as a generator of $C([0,1])$, and let
$t_1,t_2$ be generators of $\mathcal{O}_2$. Put
\[
a_1:=t_1 c,\quad a_2:=t_2(1-c).
\]
It is easy to verify that the following relations are satisfied:
\begin{equation}
\label{auxrel}
a_1^*a_1+a_2^*a_2=1,\quad a_2^*a_1=0.
\end{equation}
Then we define the images of the generators of $\mathcal{A}_2^q$ as
follows:
\begin{equation}\label{hom}
\varphi(s_1)=1\otimes s,\quad \varphi(s_2)=a_1\otimes u +a_2\otimes s
u.
\end{equation}
To verify that $\varphi(s_i^*s_i)=1$ and
$\varphi(s_2s_1)=q\varphi(s_1s_2)$, one has only to use the
relations (\ref{auxrel}).

Finally we fix the irreducible representation $\widetilde{\pi}$ of
$\mathcal{B}$ acting on $\mathcal{K}=l_2(\mathbb{N})$ by
\begin{equation}\label{rep}
\widetilde{\pi}(s)=S,\quad \widetilde{\pi}(u)=D(q),
\end{equation}
where
\[
Se_n=e_{n+1},\quad D(q)e_n=q^{n-1}e_{n},\; n\in\mathbb{N}.
\]
Then the induced functor
\[
\mathcal{F}_{\varphi}\colon \Rep \mathcal{O}_2\star
C([0,1])\rightarrow \Rep \mathcal{A}_2^q
\]
is given by the following construction.
Starting with a representation of
$\mathcal{O}_2\star C([0,1])$, say $\pi$, if as above we put
$C:=\pi(c)$, $T_i:=\pi(t_i)$, and $A_i:=\pi(a_i)$, $i=1,2$, and then we define
$\mathcal{F}_{\varphi}(\pi):=\pi_1$ by the formulas
\begin{equation}
\label{constr}
\pi_1(s_1):=S_1=1\otimes S,\quad \pi_1(s_2):=S_2=A_1\otimes
D(q)+A_2\otimes S D(q).
\end{equation}
The proof of the fullness and faithfulness of
$\mathcal{F}_{\varphi}$ is essentially the same as in the
proposition above. We note only that the equalities $C^2=A_1^*A_1$
and $T_i=A_iC^{-1}$, $i=1,2$, imply that
\[
\{A_i,A_i^*,\; i=1,2\}^{\prime}=\{T_i,T_i^*,C,\;
i=1,2\}^{\prime}=\pi(\mathcal{O}_2 \star C([0,1]))^{\prime},
\]
where again the prime denotes the commutant.
So one has to show that any
$\Lambda\in\pi_1(\mathcal{A}_2^q)^{\prime}$ has the form
$\Lambda=\Lambda_1\otimes\mathbf{1}$ with
$\Lambda_1\in\{A_i,A_i^*,\; i=1,2\}^{\prime}$.
\end{proof}
The following corollary is immediate.
\begin{corollary}
The $C^*$-algebra $\mathcal{A}_2^q$ is not nuclear.
\end{corollary}
\noindent {\bf Acknowledgements.}\\
D. Proskurin and Yu.\ Samo\u\i{}lenko express their gratitude to
Dr.\ Stanislaw Popovych for helpful discussions and consultations
during the preparation of this paper.
The authors thank Brian Treadway for suggestions, corrections, and typesetting.

\nopagebreak \vskip 0.5cm \noindent
\address
\nopagebreak

\bigskip\noindent{\em e-mail: }\email
\end{document}